\documentclass[12pt]{article}

\usepackage{amsmath,amssymb,url}

\usepackage{pict2e,graphicx}

\newtheorem{theorem}{Theorem}
\newtheorem{corollary}[theorem]{Corollary}
\newtheorem{lemma}[theorem]{Lemma}
\newtheorem{conjecture}{Conjecture}

\begin{document}

\title{Honeycomb arrays}

\author{Simon R.\ Blackburn\thanks{Department of Mathematics, Royal
Holloway, University of London, Egham, Surrey TW20 0EX, United
Kingdom. \texttt{\{s.blackburn, a.panoui\}@rhul.ac.uk}}\\
Anastasia Panoui$^*$\\
Maura B. Paterson\thanks{Department of Economics, Mathematics and Statistics, Birkbeck, University of London, Malet Street, London WC1E 7HX, United Kingdom. \texttt{m.paterson@bbk.ac.uk}}\\ and\\
Douglas R. Stinson\thanks{David R.\ Cheriton School of Computer Science, University of Waterloo, Waterloo Ontario, N2L 3G1, Canada. \texttt{dstinson@uwaterloo.ca}}}

\maketitle

\begin{abstract}
A honeycomb array is an analogue of a Costas array in the hexagonal
grid; they were first studied by Golomb and Taylor in 1984. A recent
result of Blackburn, Etzion, Martin and Paterson has shown that (in
contrast to the situation for Costas arrays) there are only finitely
many examples of honeycomb arrays, though their bound on the maximal
size of a honeycomb array is too large to permit an exhaustive search
over all possibilities.

The present paper contains a theorem that significantly limits the
number of possibilities for a honeycomb array (in particular, the
theorem implies that the number of dots in a honeycomb array must be
odd). Computer searches for honeycomb arrays are summarised, and two
new examples of honeycomb arrays with 15 dots are given.
\end{abstract}

\section{Introduction}

Honeycomb arrays were introduced by Golomb and
Taylor~\cite{GolombTaylor} in 1984, as a hexagonal analogue of Costas
arrays. Examples of honeycomb arrays are given in
Figures~\ref{honey137_figure} to~\ref{honey45_figure} below.
\begin{figure}
\begin{center}
\includegraphics[width=70mm, trim=0mm 25mm 0mm 25mm]{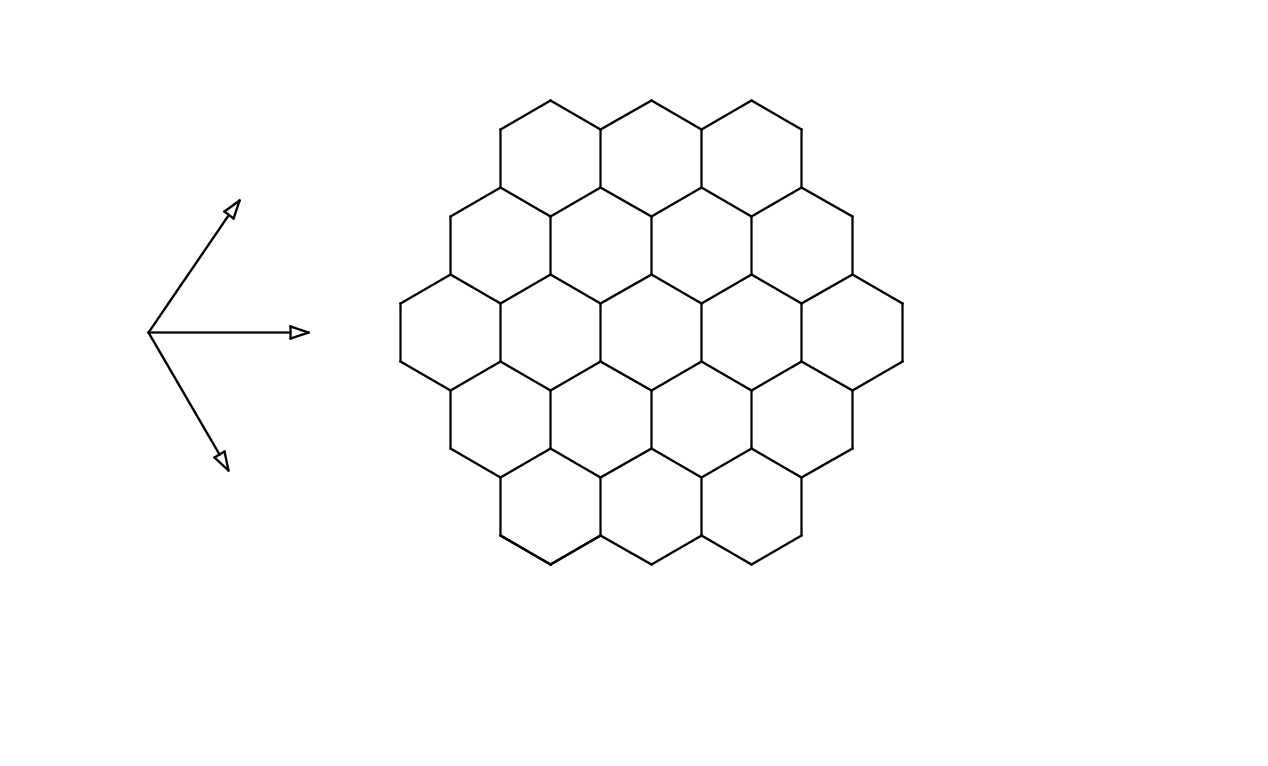}
\end{center}
\caption{A Lee sphere, and three natural directions}
\label{fig:Lee_example}
\end{figure}A honeycomb array is a
collection of $n$ dots in a hexagonal array with two properties:
\begin{itemize}
\item \textbf{(The hexagonal permutation property)} There are three natural
directions in a hexagonal grid (see
Figure~\ref{fig:Lee_example}). Considering `rows' in each of these
three directions, the dots occupy $n$ consecutive rows, with exactly
one dot in each row.
\item \textbf{(The distinct differences property)} The $n(n-1)$ vector
differences between pairs of distinct dots are all different.
\end{itemize}

Golomb and Taylor found $10$ examples of honeycomb arrays (up to
symmetry), and conjectured that infinite families of honeycomb arrays
exist. Blackburn, Etzion, Martin and Paterson~\cite{BlackburnEtzion}
recently disproved this conjecture: there are only a finite number of
honeycomb arrays. Unfortunately, the bound on the maximal size of a
honeycomb array that Blackburn et al.\ provide is far too large to
enable an exhaustive computer search over all open cases. In this
paper, we prove a theorem that significantly limits the possibilities
for a honeycomb array with $n$ dots. (In particular, we show that $n$
must be odd.) We report on our computer searches for honeycomb arrays,
and give two previously unknown examples with $15$ dots.

We now introduce a little more notation, so that we can state the main
result of our paper more precisely.

We say that a collection of dots in the hexagonal grid is a
\emph{hexagonal permutation} if it satisfies the hexagonal permutation
property. A collection of dots is a \emph{distinct difference
configuration} if it satisfies the distinct difference property. So a
honeycomb array is a hexagonal permutation that is a distinct
difference configuration.

We say that hexagons are \emph{adjacent} if they share an edge, and we
say that two hexagons $A$ and $B$ are \emph{at distance $d$} if the
shortest path from $A$ to $B$ (travelling along adjacent hexagons) has
length $d$. A \emph{Lee sphere of radius $r$} is a region of the
hexagonal grid consisting of all hexagons at distance $r$ or less from
a fixed hexagon (the \emph{centre} of the sphere). The region in
Figure~\ref{fig:Lee_example} is a Lee sphere of radius $2$. Note that
a Lee sphere of radius $r$ intersects exactly $2r+1$ rows in each of
the three natural directions in the grid. A \emph{honeycomb array of
radius $r$} is a honeycomb array with $2r+1$ dots contained in a Lee
sphere of radius $r$.

\begin{figure}
\begin{center}
\includegraphics[width=70mm,trim=0mm 20mm 0mm 20mm]{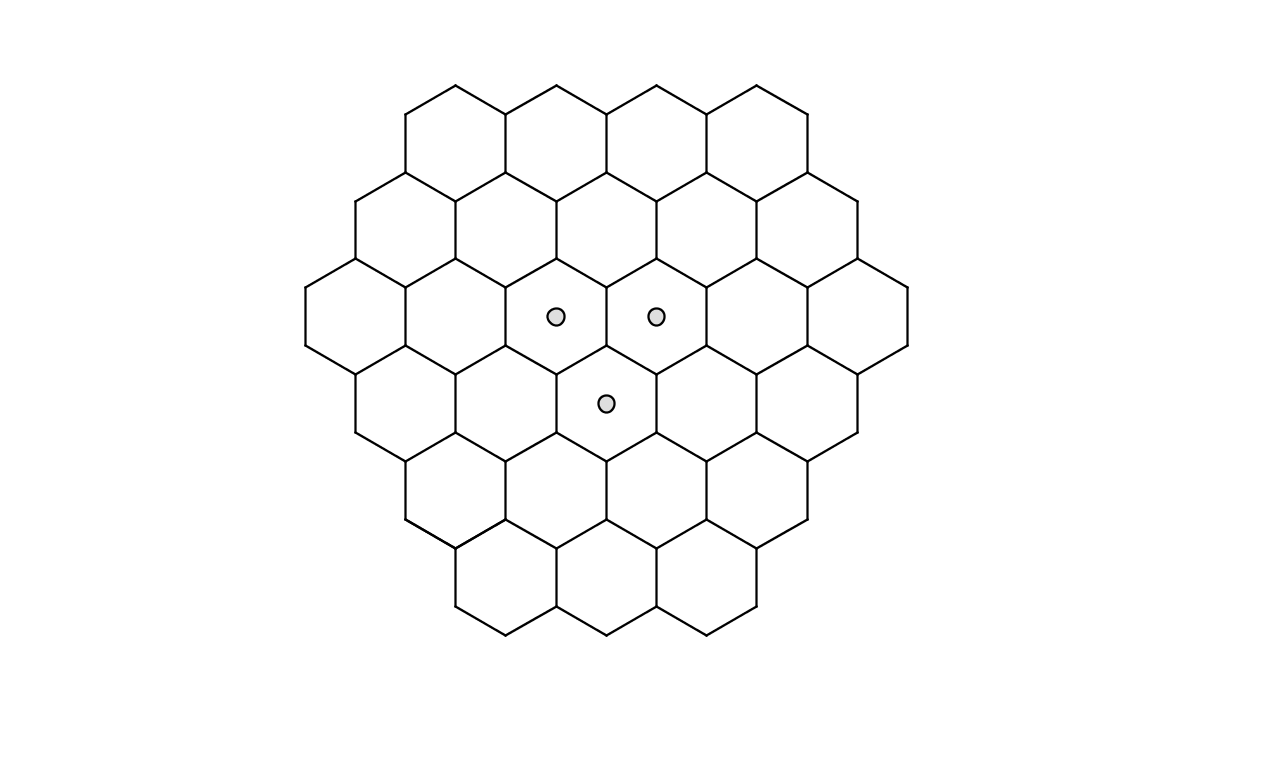}
\end{center}
\caption{A tricentred Lee sphere}
\label{fig:tricentred_example}
\end{figure}

There are many other natural regions of the hexagonal grid that have
the property that they intersect $n$ rows in each direction. One
example, the tricentred Lee sphere of radius $r$, is shown in
Figure~\ref{fig:tricentred_example}: it is the union of three Lee
spheres of radius $r$ with pairwise adjacent centres, and intersects
exactly $2r+2$ rows in any direction.

Does there exist a honeycomb
array with $2r+2$ dots contained in a tricentred Lee sphere of
radius~$r$? Golomb and Taylor did not find any such examples: they
commented~\cite[Page~1156]{GolombTaylor} that all known examples of
honeycomb arrays with $n$ dots were in fact honeycomb arrays of radius
$r$, but stated ``we have not proved that this must always be the
case''. We prove the following:

\begin{theorem}
\label{thm:main}
Let $n$ be an integer, and suppose there exists a hexagonal
permutation $\pi$ with $n$ dots. Then $n$ is odd, and the dots of $\pi$ are
contained in a Lee sphere of radius $(n-1)/2$.
\end{theorem}

Since any honeycomb array is a hexagonal permutation, the following
result follows immediately from Theorem~\ref{thm:main}:

\begin{corollary}
\label{cor:honeycomb}
Any honeycomb array is a honeycomb array of radius $r$ for some
integer $r$. In particular, a honeycomb array must consist of an odd
number of dots.
\end{corollary}

So if we are looking for honeycomb arrays, we may restrict ourselves to
searching for honeycomb arrays of radius $r$.

The structure of the remainder of the paper is as follows. In
Section~\ref{sec:anticodes}, we state the results on the hexagonal
grid that we need. In Section~\ref{sec:brooks}, we remind the reader
of the notion of a brook, or bee-rook, and state a theorem on the
maximum number of non-attacking brooks on a triangular board. In
Section~\ref{sec:honeycomb} we prove Theorem~\ref{thm:main}, we
summarise our computer searches for honeycomb arrays, and we provide a
list of all known arrays. This last section also contains a conjecture, and
some suggestions for further work.

\section{The hexagonal grid}
\label{sec:anticodes}

\begin{figure}
\centering \setlength{\unitlength}{.55mm}
\begin{picture}(150,55)(-10,-40)
\linethickness{.1pt} \multiput(10,5)(10,0){3}{
\put(-10,0){\line(0,1){5.77350269}}
\put(0,0){\line(0,1){5.77350269}}
\put(0,5.77350269){\line(-1732,1000){5}}
\put(-10,5.77350269){\line(1732,1000){5}} }
\multiput(5,-3.66033872)(10,0){4}{
\put(-10,0){\line(0,1){5.77350269}}
\put(0,0){\line(0,1){5.77350269}}
\put(0,5.77350269){\line(-1732,1000){5}}
\put(-10,5.77350269){\line(1732,1000){5}} }
\multiput(0,-12.3206774)(10,0){5}{
\put(-10,0){\line(0,1){5.77350269}}
\put(0,0){\line(0,1){5.77350269}}
\put(0,5.77350269){\line(-1732,1000){5}}
\put(-10,5.77350269){\line(1732,1000){5}}
}
\multiput(5, -20.9810162)(10,0){4}{
\put(-10,0){\line(0,1){5.77350269}}
\put(0,0){\line(0,1){5.77350269}}
\put(0,5.77350269){\line(-1732,1000){5}}
\put(-10,5.77350269){\line(-1732,1000){5}}
\put(-10,5.77350269){\line(1732,1000){5}}
\put(0,5.77350269){\line(1732,1000){5}}
}
\multiput(10, -29.6413549)(10,0){3}{
\put(-10,0){\line(0,1){5.77350269}}
\put(0,0){\line(0,1){5.77350269}}
\put(0,5.77350269){\line(-1732,1000){5}}
\put(-10,5.77350269){\line(-1732,1000){5}}
\put(-10,5.77350269){\line(1732,1000){5}}
\put(0,5.77350269){\line(1732,1000){5}}
}
\multiput(15, -38.3016936)(10,0){2}{
\put(0,5.77350269){\line(-1732,1000){5}}
\put(-10,5.77350269){\line(-1732,1000){5}}
\put(-10,5.77350269){\line(1732,1000){5}}
\put(0,5.77350269){\line(1732,1000){5}}}

\put(15,-9.43392605){\makebox(0,0){$0$}}
\put(5,-9.43392605){\makebox(0,0){$5$}}
\put(25,-9.43392605){\makebox(0,0){$2$}}
\put(10,-0.77358733){\makebox(0,0){$6$}}
\put(20,-0.77358733){\makebox(0,0){$1$}}
\put(10,-18.0942648){\makebox(0,0){$4$}}
\put(20,-18.0942648){\makebox(0,0){$3$}}
%
\thinlines


\put(60,-16.43392605){\makebox(0,0){$\overrightarrow{\quad\xi\quad}$}}

\put(80,-38.3014395){\line(1,0){34.6410162}}
\put(80,-26.7544341){\line(1,0){46.1880215}}
\put(80,-15.2074287){\line(1,0){57.7350269}}
\put(80.5470054,-3.66042332){\line(1,0){57.7350269}}
\put(91.5470054,7.88658206){\line(1,0){46.1880215}}
\put(103.094011, 19.4335874){\line(1,0){34.6410162}}

\put(80,-38.3014395){\line(0,1){34.6410162}}
\put(91.5470054,-38.3014395){\line(0,1){46.1880215}}
\put(103.094011,-38.3014395){\line(0,1){57.7350269}}
\put(114.641016,-38.3014395){\line(0,1){57.7350269}}
\put(126.188022,-26.7544341){\line(0,1){46.1880215}}
\put(137.735027,-15.2074287){\line(0,1){34.6410162}}

\put(108.867513,-9.43392605){\makebox(0,0){$0$}}
\put( 97.3205081,-9.43392605){\makebox(0,0){$5$}}
\put(120.414519,-9.43392605){\makebox(0,0){$2$}}
\put(108.867513,2.11307933){\makebox(0,0){$6$}}
\put(97.3205081,-20.9809314){\makebox(0,0){$4$}}
\put(120.414519,2.11307933){\makebox(0,0){$1$}}
\put(108.867513,-20.9809314){\makebox(0,0){$3$}}
\end{picture}
\caption{From the hexagonal to the square grid}
\label{fig:hex_to_square}
\end{figure}
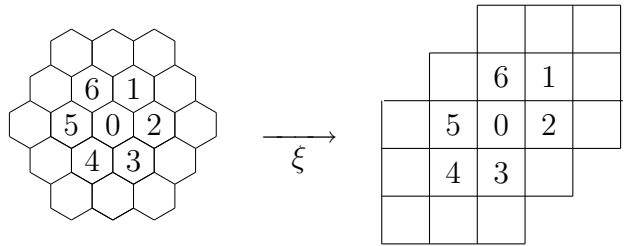

Because the hexagonal grid might be difficult to visualise, we use an
equivalent representation in the square grid (see
Figure~\ref{fig:hex_to_square}). In this representation, we define
each square to be adjacent to the four squares it shares an edge with,
and the squares sharing its `North-East' and `South-West' corner
vertices. The map $\xi$ in Figure~\ref{fig:hex_to_square} distorts the
centres of the hexagons in the hexagonal grid to the centres of the
squares in the square grid. The three types of rows in the hexagonal
grid become the rows, the columns and the diagonals that run from
North-East to South-West. For brevity, we define a \emph{standard
diagonal} to mean a diagonal that runs North-East to South-West.

\begin{figure}
\begin{center}
\includegraphics[width=70mm, trim=0mm 8mm 0mm 8mm]{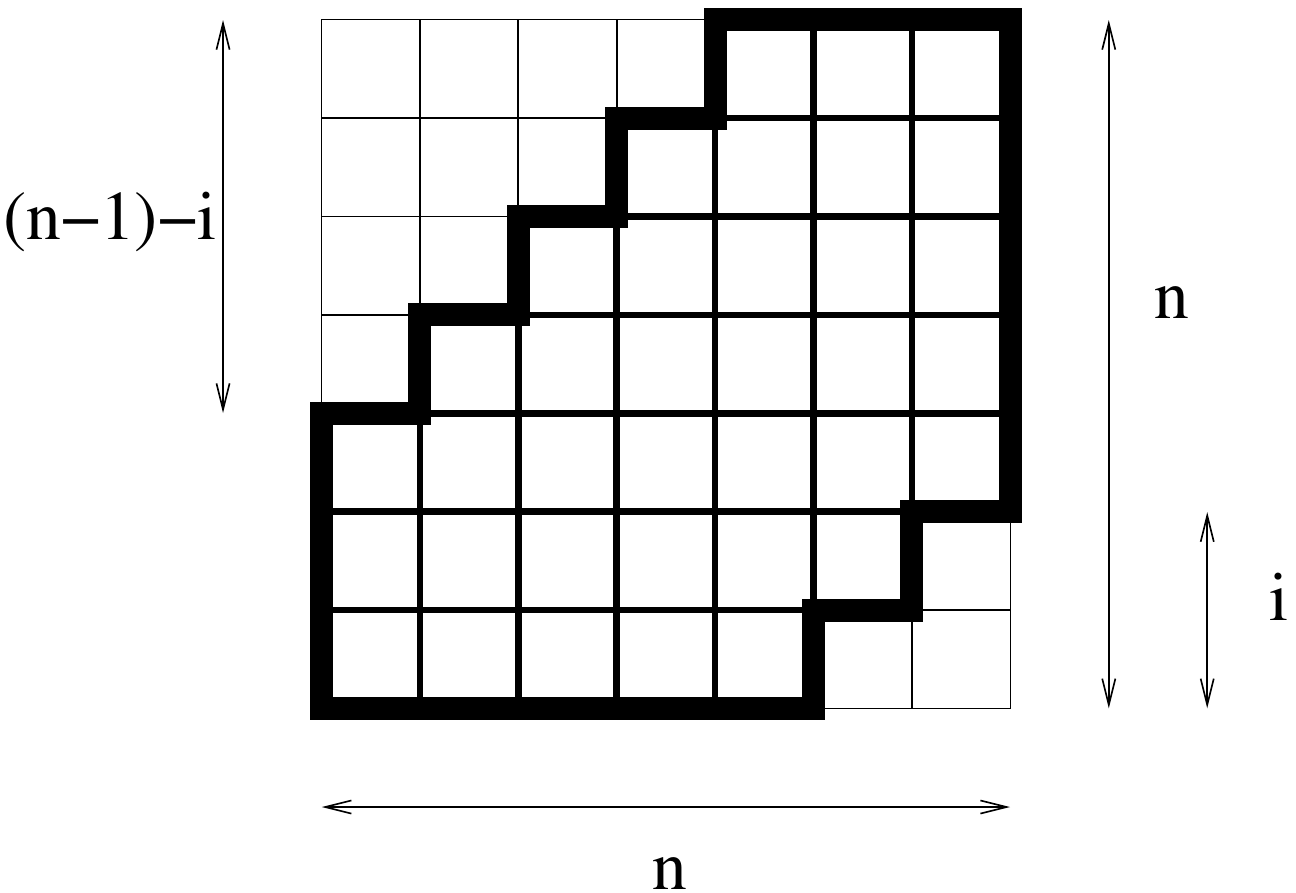}
\end{center}
\caption{The region $S_{i}(n)$}
\label{fig:anticode}
\end{figure}
For non-negative integers $n$ and $i$ such that $0\leq i\leq n-1$,
define the region $S_i(n)$ of the square grid as in Figure~\ref{fig:anticode}. Note that
$S_{i}(n)$ and $S_{n-1-i}(n)$ are essentially the same region: one is
obtained from the other by a reflection in a standard diagonal. The
regions $\xi^{-1}(S_i(n))$ are important in the hexagonal grid, as they are the maximal
anticodes of diameter $n-1$; see Blackburn et
al.~\cite[Theorem~5]{BlackburnEtzion}.  Note that the region
$\xi^{-1}(S_{r}(2r+1))$ is a Lee sphere of radius $r$. Regions of the form
$\xi^{-1}(S_{r}(2r+2))$ or $\xi^{-1}(S_{r+1}(2r+2))$ are tricentred Lee spheres of radius
$r$. Also note that the regions $S_i(n)$ as $i$ varies are precisely
the possible intersections of an $n\times n$ square region with $n$
adjacent standard diagonals, where each diagonal intersects the
$n\times n$ square non-trivially.

In the lemma below, by a `region of the form $X$', we
mean a region that is a translation of $X$ in the square grid.

\begin{lemma}
\label{lem:anticode}
Let $\pi$ be a hexagonal permutation with $n$ dots, and let $\xi(\pi)$ be the image of $\pi$ in the square grid.
Then the dots in
$\xi(\pi)$ are all contained in a region of the form $S_{i}(n)$ for some
$i$ in the range $0\leq i\leq n-1$.
\end{lemma}
\textbf{Proof:} Let $R$ be the set of squares that share a row with a
dot of $\xi(\pi)$. Similarly, let $C$ and $D$ be the sets squares sharing
respectively a column or a standard diagonal with a dot of $\xi(\pi)$. The
dots in $\xi(\pi)$ are contained in $R\cap C\cap D$.

Since $\pi$ is a hexagonal permutation, $R$ consists of $n$ adjacent
rows and $C$ consists of $n$ adjacent columns. Hence $R\cap C$ is an
$n\times n$ square region. (Since there is exactly one dot in each row
and column of the square $R\cap C$, the dots in $\xi(\pi)$ correspond to a
permutation; this justifies the terminology `hexagonal permutation'.)

Now, $D$ consists of $n$ adjacent standard diagonals; each diagonal contains a
dot in $\xi(\pi)$, and so each diagonal intersects $R\cap C$
non-trivially. Hence $R\cap C\cap D$ is a region of the form $S_i(n)$,
as required.\hfill$\Box$

\section{Brooks on a triangular board}
\label{sec:brooks}

A \emph{brook} is a chess piece in the square grid that moves like a rook plus half a
bishop: it can move any distance along a row, a column or a standard
(North-East to South-West) diagonal. Brooks were first studied by
Bennett and Potts~\cite{BennettPotts}, who pointed out connections to
constant-sum arrays and hexagonal lattices. A set of brooks in a
square grid is \emph{non-attacking} if no two brooks lie in a row, a
column or a standard diagonal.

Under the correspondence $\xi$ between the square and hexagonal grids
mentioned in the previous section, brooks in the square grid
correspond to \emph{bee-rooks} in the hexagonal grid: pieces that can move
any distance along any row, where a row can go in each of the three
natural directions. A set of bee-rooks is therefore \emph{non-attacking}
if no two bee-rooks lie in the same row of the hexagonal grid. In particular,
bee-rooks placed on the dots in a hexagonal permutation $\pi$
are non-attacking, and so the corresponding set $\xi(\pi)$ of brooks in the
square grid is non-attacking.

A \emph{triangular board of width $w$} is the region $S_{0}(w)$ in the square grid
depicted in Figure~\ref{fig:triangular_board}. Let $b(w)$ be the
maximum number of non-attacking brooks that can be placed in the
triangular board of width $w$. The following theorem is proved by
Nivasch and Lev~\cite{NivaschLev} and in Vaderlind, Guy and
Larson~\cite[P252 and R252]{VaderlindGuy}:
\begin{figure}
\begin{center}
\includegraphics[width=50mm, trim=0mm 8mm 0mm 8mm]{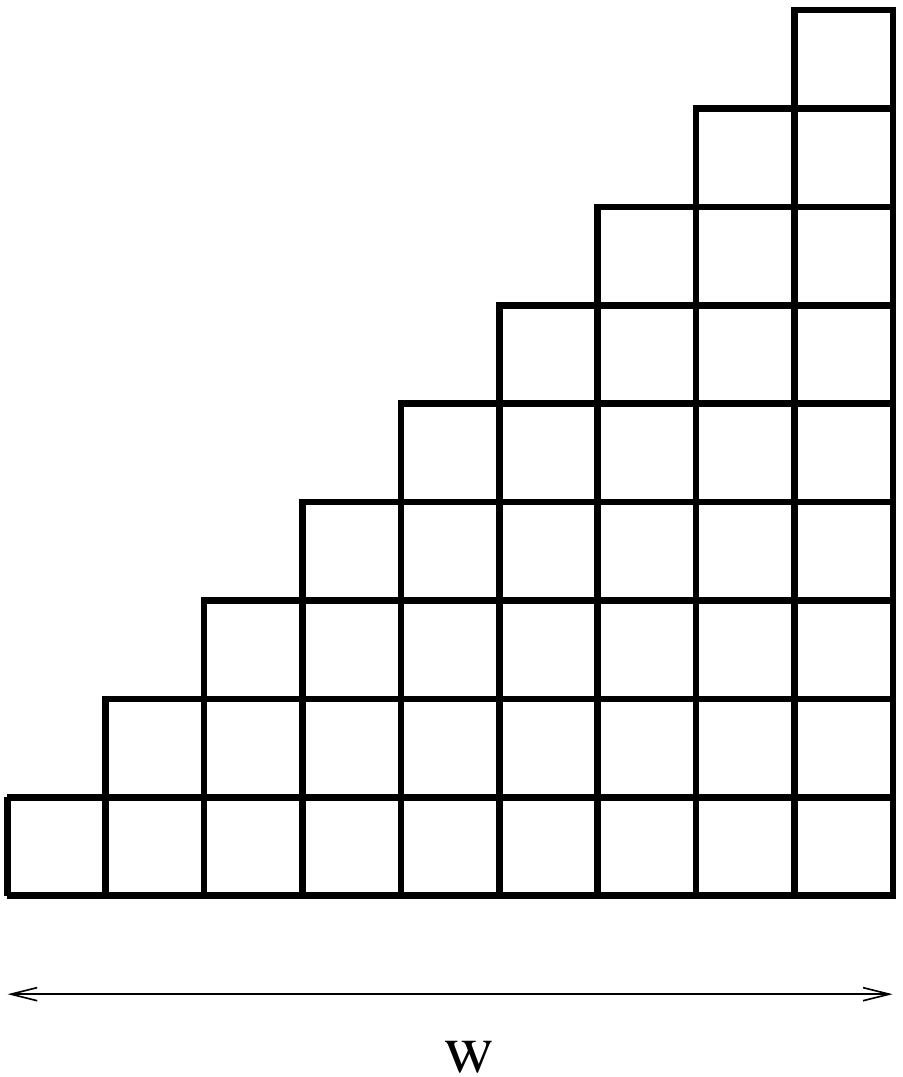}
\end{center}
\caption{A triangular board of width $w$}
\label{fig:triangular_board}
\end{figure}

\begin{theorem}
\label{thm:brooks}
For any positive integer $w$, $b(w)=\lfloor (2w+1)/3\rfloor$.
\end{theorem}

Three of the present authors have found an alternative proof for this
theorem, using linear programming
techniques~\cite{BlackburnPaterson}. See Bell and
Stevens~\cite{BellStevens} for a survey of similar combinatorial
problems.

\section{Honeycomb arrays}
\label{sec:honeycomb}

We begin this section with a proof of Theorem~\ref{thm:main}. We then
describe our searches for honeycomb arrays. We end the section by
describing some avenues for further work.

\noindent
\textbf{Proof of Theorem~\ref{thm:main}:} Let $\pi$ be a hexagonal
permutation with $n$ dots. By
Lemma~\ref{lem:anticode}, the dots of $\xi(\pi)$ are contained in a region
of the form $S_i(n)$ where $0\leq i\leq n-1$. When $i=(n-1)/2$ (so $n$
is odd and $\xi^{-1}(S_i(n))$ is a Lee sphere) the
theorem follows. Suppose, for a contradiction, that
$i\not=(n-1)/2$.

By reflecting $\pi$ in a horizontal row in the hexagonal grid, we produce a hexagonal
permutation $\pi'$ such that $\xi(\pi')$ is contained in a region of the form
$S_{(n-1)-i}(n)$. By replacing $\pi$ by $\pi'$ if necessary, we may
assume that $i<(n-1)/2$.

\begin{figure}
\begin{center}
\includegraphics[width=65mm, trim=0mm 8mm 0mm 8mm]{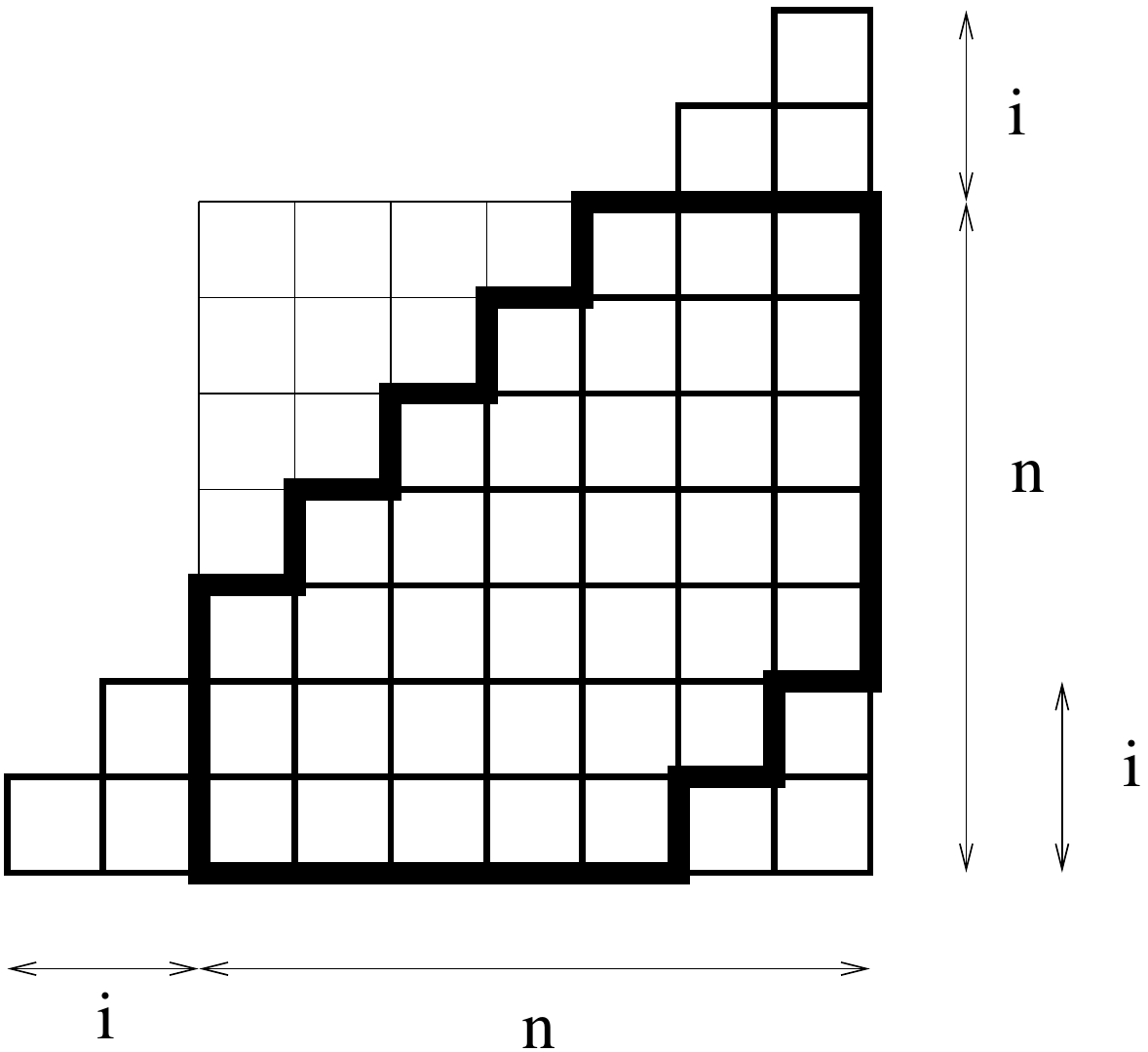}
\end{center}
\caption{A triangular board covering $\pi$}
\label{fig:covering_triangle}
\end{figure}
Consider the triangular board of width $n+i$ in
Figure~\ref{fig:covering_triangle} containing $S_i(n)$. Since no two
dots in $\xi(\pi)$ lie in the same row, column or standard diagonal, the
dots in $\xi(\pi)$ correspond to $n$ non-attacking brooks in this
triangular board. But this contradicts Theorem~\ref{thm:brooks}, since
\[
\frac{2(n+i)+1}{3}<\frac{2n+(n-1)+1}{3}= n.
\]
This contradiction completes the proof of the theorem.\hfill$\Box$

\begin{figure}
\begin{center}
\parbox[c]{10mm}{\includegraphics[width=6mm,angle=90]{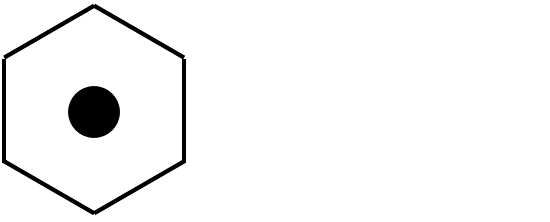}}
\parbox[c]{25mm}{\includegraphics[width=17mm,angle=90]{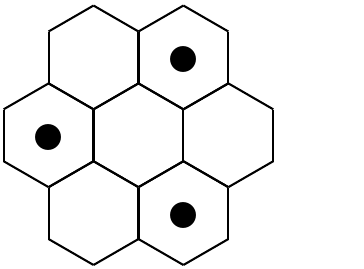}}
\parbox[c]{45mm}{\includegraphics[width=35mm,angle=90]{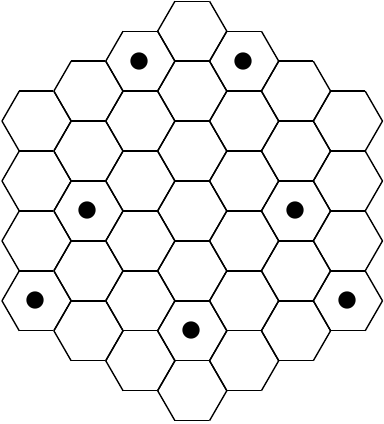}}
\parbox[c]{45mm}{\includegraphics[width=35mm,angle=90]{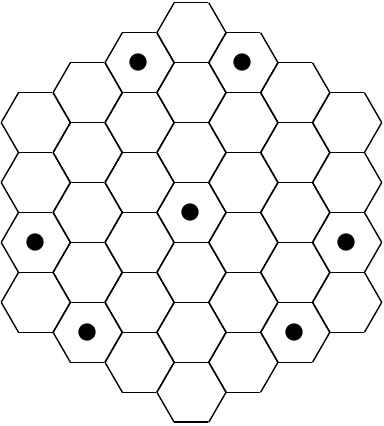}}
\end{center}
\caption{Honeycomb arrays of radius $0$, $1$ and $3$}
\label{honey137_figure}
\end{figure}

\begin{figure}
\begin{center}
\parbox[c]{50mm}{\includegraphics[width=40mm,angle=90]{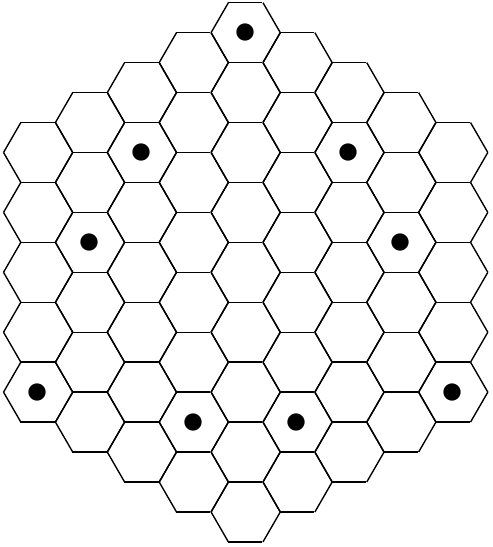}}
\parbox[c]{50mm}{\includegraphics[width=40mm,angle=90]{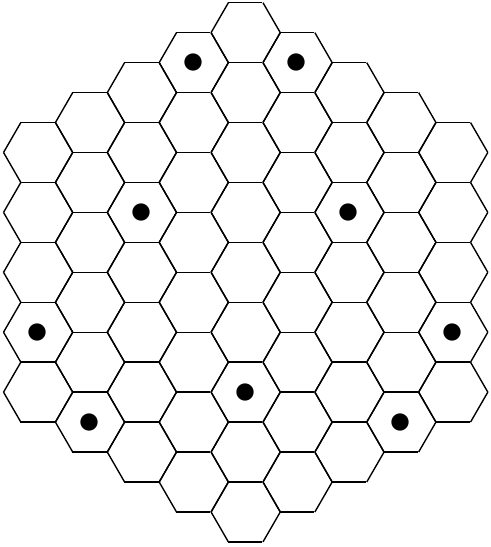}}
\end{center}
\caption{Honeycomb arrays of radius 4}
\label{honey9_figure}
\end{figure}

\begin{figure}
\begin{center}
\parbox[c]{40mm}{\includegraphics[width=35mm,angle=90]{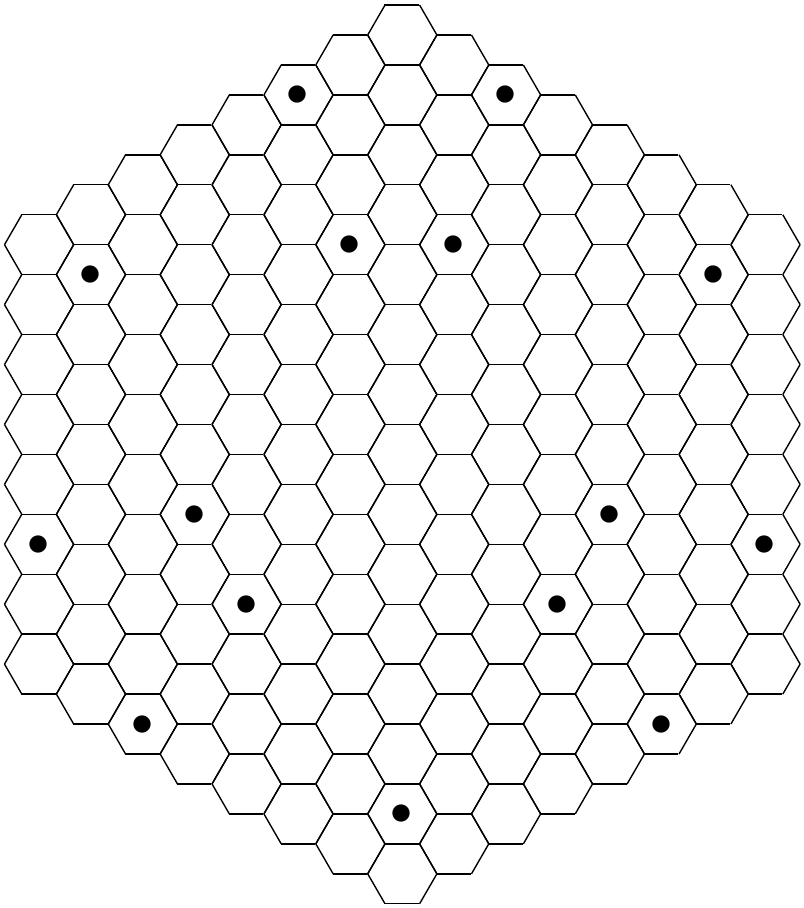}}
\parbox[c]{40mm}{\includegraphics[width=35mm,angle=90]{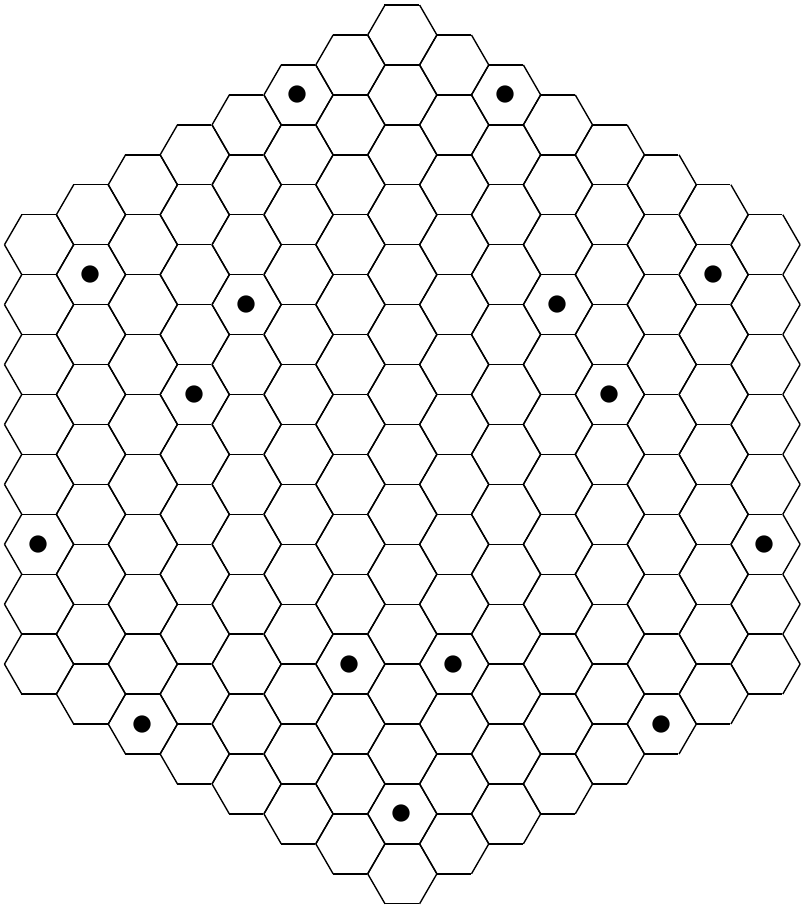}}
\parbox[c]{40mm}{\includegraphics[width=35mm,angle=90]{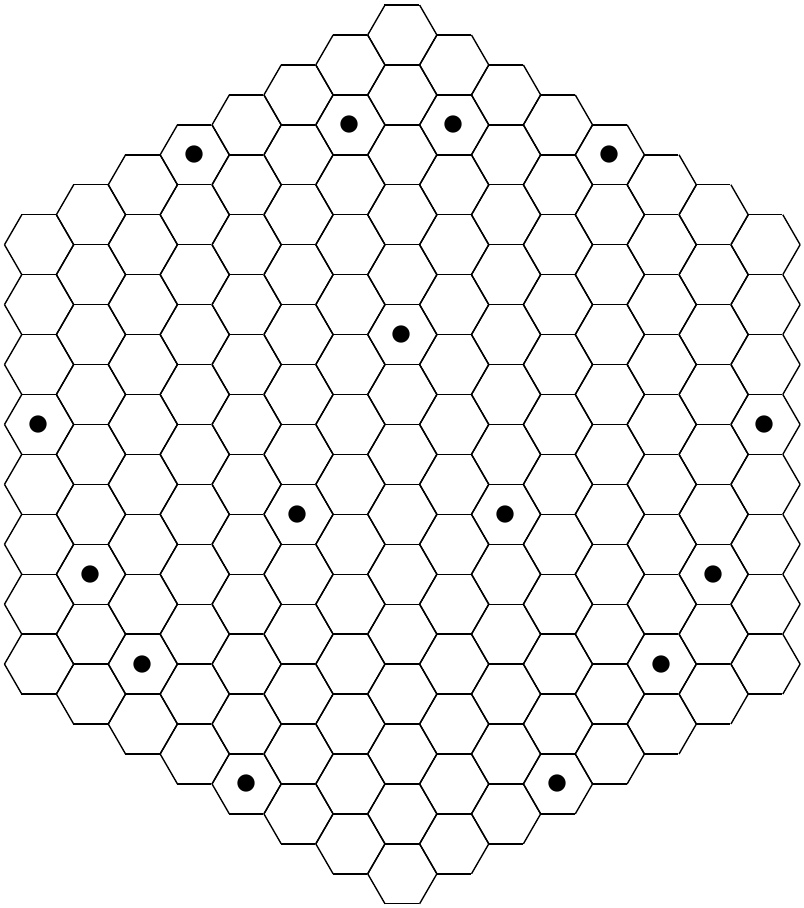}}
\end{center}
\caption{Honeycomb arrays of radius $7$}
\label{honey15_figure}
\end{figure}

\begin{figure}
\begin{center}
\parbox[c]{60mm}{\includegraphics[height=58mm,angle=90]{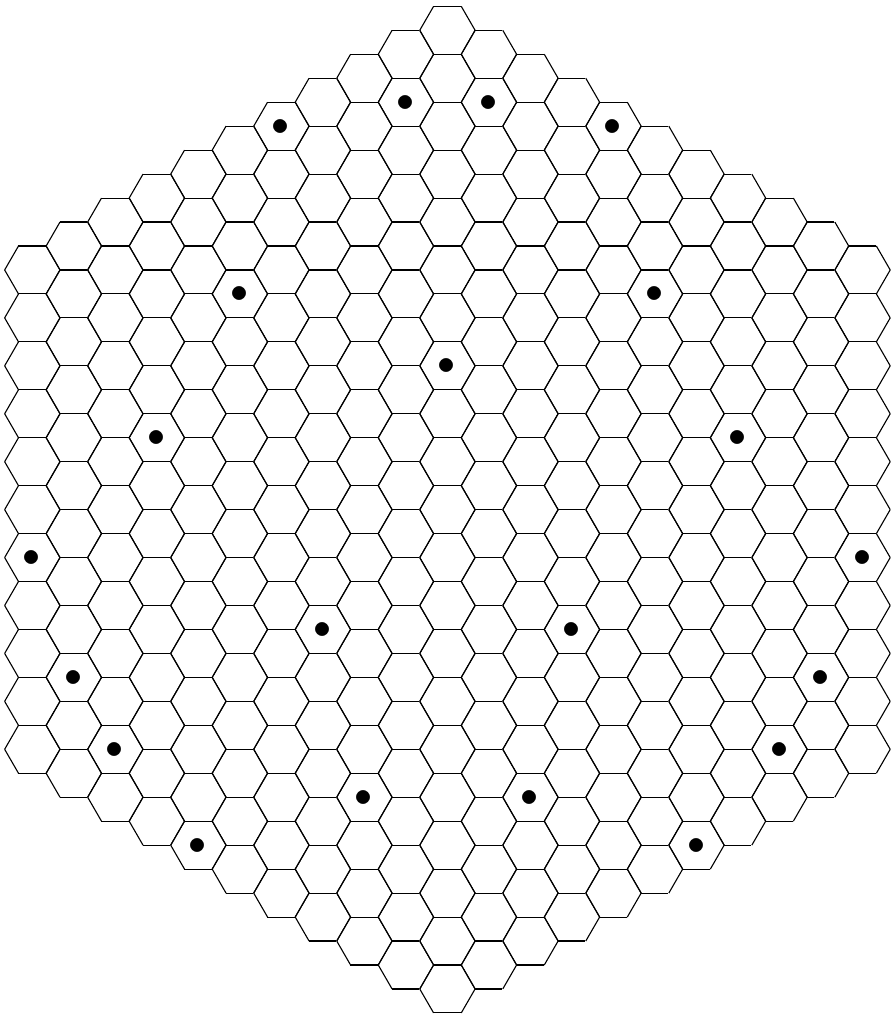}}
\parbox[c]{75mm}{\includegraphics[height=73mm,angle=90]{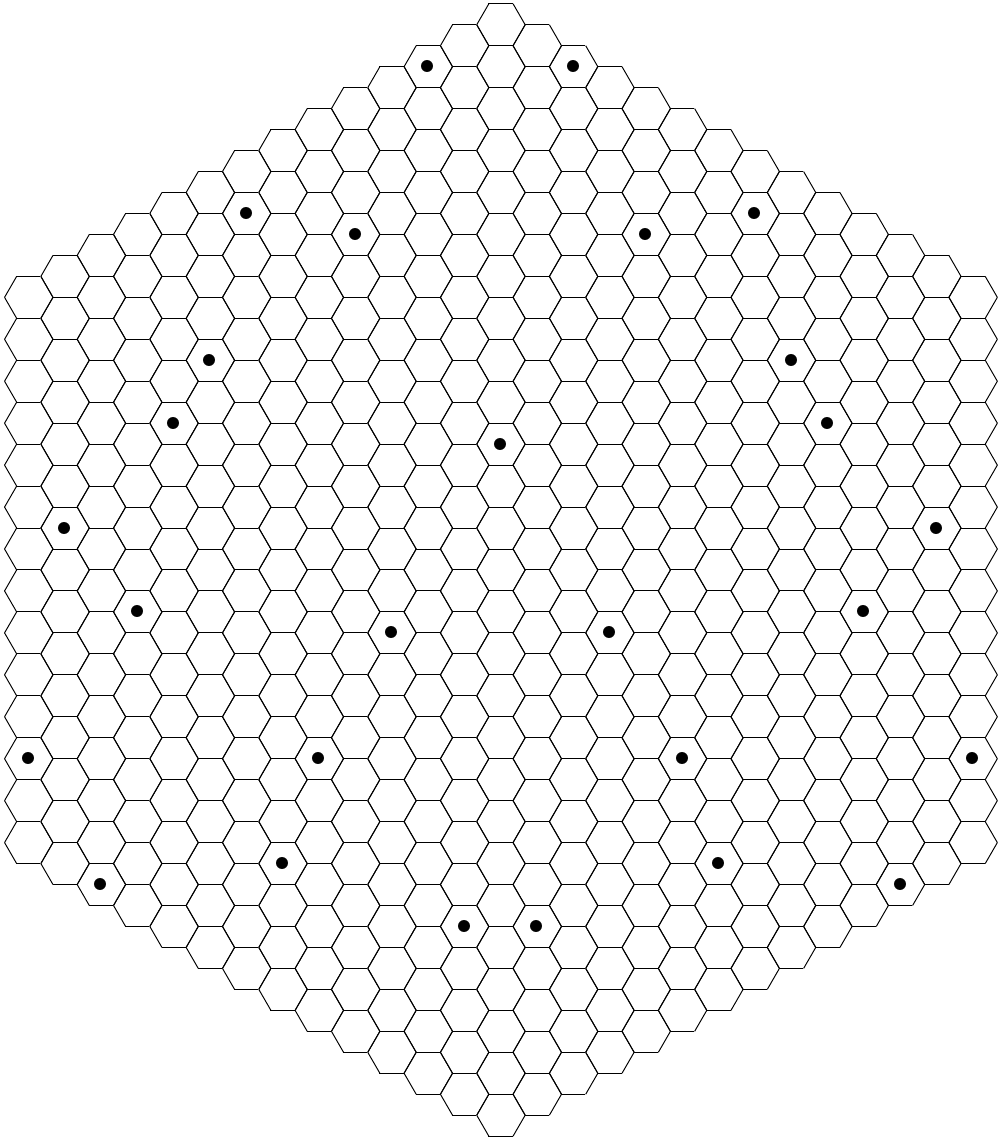}}
\end{center}
\caption{Honeycomb arrays of radius $10$ and $13$}
\label{honey2127_figure}
\end{figure}

\begin{figure}
\begin{center}
\includegraphics[height=140mm,angle=90]{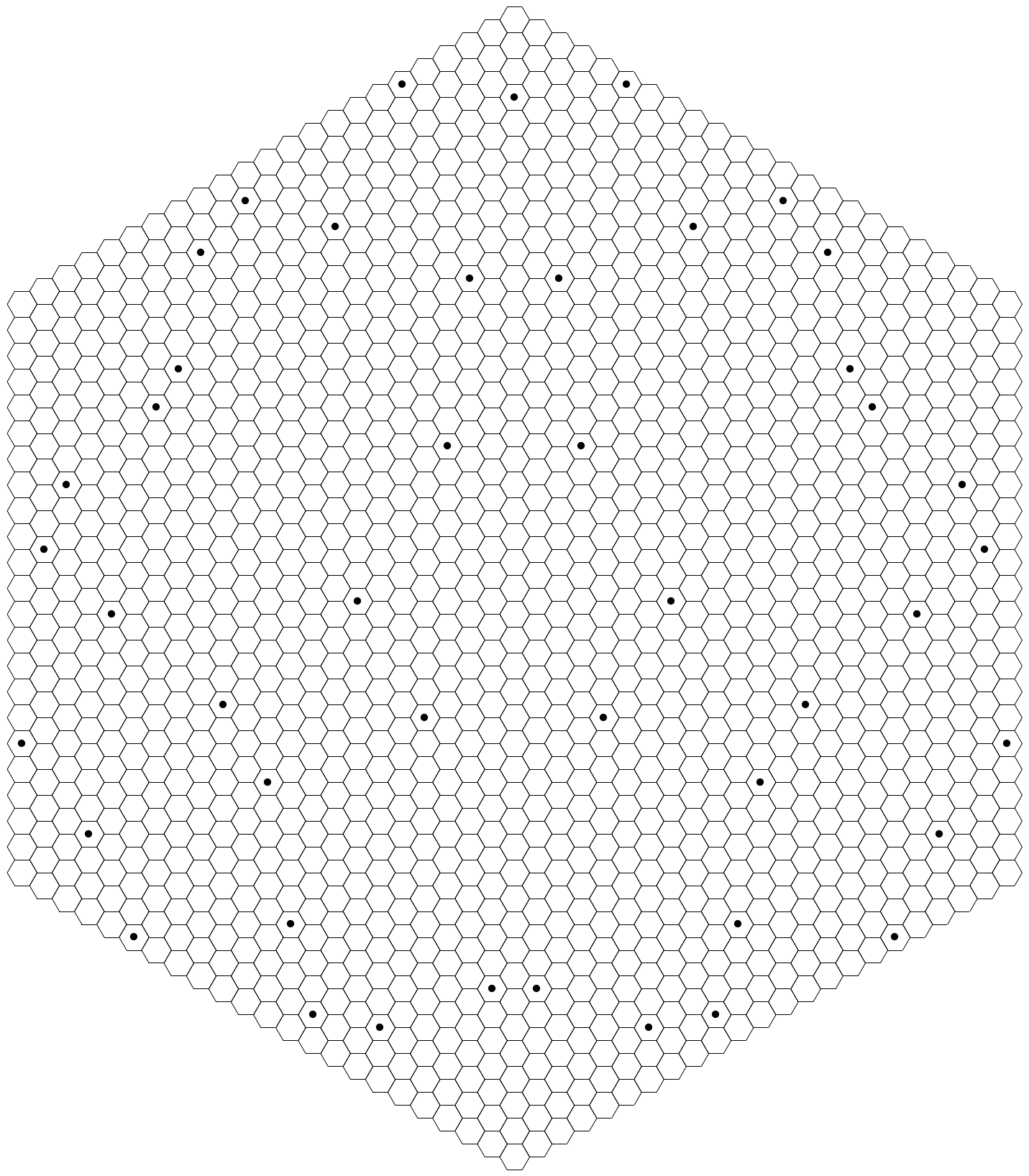}
\end{center}
\caption{A honeycomb array of radius $22$}
\label{honey45_figure}
\end{figure}

Theorem~\ref{thm:main} tells us that the only honeycomb arrays are
those of radius $r$ for some non-negative integer $r$. A result of
Blackburn et al~\cite[Corollary~12]{BlackburnEtzion} shows that $r\leq
643$. We now report on our computer searches for examples of honeycomb
arrays. The known honeycomb arrays
are drawn in Figures~\ref{honey137_figure}, \ref{honey9_figure},
\ref{honey15_figure}, \ref{honey2127_figure} and~\ref{honey45_figure}.
This list includes two new examples not known to Golomb and
Taylor~\cite{GolombTaylor}, namely the second and third examples of
radius $7$; we found these examples as follows.

A \emph{Costas array} is a set of $n$ dots in an $n\times n$ region of
the square grid, with the distict difference property and such that
every row and column of the array contains exactly one dot. Golomb and
Taylor observed that some Costas arrays produce honeycomb arrays, by
mapping the dots in the Costas array into the hexagonal grid using the
map $\xi^{-1}$ given by Figure~\ref{fig:hex_to_square}. Indeed, it is
not difficult to see that all honeycomb arrays must arise in this
way. We searched for honeycomb arrays by taking each known Costas
array with $200$ or fewer dots, and checking whether the array gives
rise to a honeycomb array. For our search, we made use of a database
of all known Costas arrays with 200 or fewer dots that has been
made available by James K. Beard~\cite{Beard}. This list is known to be
complete for Costas arrays with $27$ or fewer dots; see Drakakis et
al.~\cite{DrakakisRickard} for details. So our list of honeycomb
arrays of radius $13$ or less is complete.

It is a remarkable fact that all known honeycomb arrays possess a
non-trivial symmetry (a horizontal reflection as we have drawn
them). Indeed, apart from a single example of radius $3$ (the
first radius $3$ example in Figure~\ref{honey137_figure}) all
known honeycomb arrays possess a symmetry group of order $6$: the group
generated by the reflections along the three lines through
opposite `corners' of the hexagonal sphere. We implemented an
exhaustive search for honeycomb arrays with $r\leq 31$ having this
$6$-fold symmetry: we found no new examples. We also checked all
constructions of honeycomb arrays from Costas arrays in Golomb and
Taylor~\cite{GolombTaylor} (whether symmetrical or not) for $r\leq
325$, and again we found no new examples.

After these searches, we feel that we can make the following conjecture:

\begin{conjecture}
The list of known honeycomb arrays is complete. So there are are
exactly $12$ honeycomb arrays, up to symmetry.
\end{conjecture}

Theorem~\ref{thm:main} shows that hexagonal permutations are always
contained in some Lee sphere. But such permutations have been
prevously studied in several contexts: Bennett and
Potts~\cite{BennettPotts} study them as non-attacking configurations
of bee-rooks and as the number of zero-sum arrays; Kotzig and
Laufer~\cite{Kotzig} study them as the number of
$\sigma$-permutations; Bebeacua, Mansour, Postnikov and
Severini~\cite{BebeacuaMansour} study them as X-rays of permutations
with maximum degeneracy. Let $h_n$ be the number of hexagonal
permutations with $2n-1$ dots. The On-Line Encyclopedia of Integer
Sequences~\cite[Sequence A002047]{OEIS} quotes a computation due to
Alex Fink that computes the first few terms of the sequence $h_n$:
\[
\begin{array}{c|cccccccccc}
n&1&2&3&4&5&6&7&8&9&10\\\hline
h_n&1&2&6&28&244&2544&35600&659632&15106128&425802176
\end{array}
\]
Kotzig and Laufer ask: How big can $h_n$ be? It seems that the sequence
grows faster than exponentially with $n$. We ask a more precise
question: Is it true that $(\log h_n)/n\log n$ tends to a constant as
$n\rightarrow\infty$?

\paragraph{Acknowledgements} Part of this work was completed under
EPSRC Grant EP/D053285/1. The authors would like to thank Tuvi Etzion for
discussions, funded by a Royal Society International Travel Grant, which
inspired this line of research.

\end{document}